\documentclass[12pt]{article}
\usepackage{amssymb}
\usepackage{theorem}
\theorembodyfont{\rmfamily}
\theoremheaderfont{\rmfamily}
\newtheorem{rem}{\rmfamily Remark}[section]
\newtheorem{dfn}{\rmfamily Definition}[section]
\newtheorem{ex}{\rmfamily Example}[section]
\newtheorem{thm}{\bfseries Theorem}[section]
\newtheorem{prop}{\bfseries Proposition}[section]
\newtheorem{lem}{\bfseries Lemma}[section]
\newtheorem{cor}{\bfseries Corollary}[section]
\theorembodyfont{\itshape}

\def\A{{\cal{A}}}
\def\B{{\cal{B}}}
\def\E{{\cal{E}}}
\def\te{{\widetilde{\cal{E}}}}
\def\tb{{\widetilde{\cal{B}}}}

\def\R{{\bf{R}}}
\def\Z{{\bf{Z}}}

\def\h{{\mathfrak{h}}}

\def\pa{{\partial}}
\def\S{{\mathfrak{S}}}
\def\SS{{\bf{S}}}
\def\lla{{\langle \! \langle}}
\def\rra{{\rangle \! \rangle}}
\title{Extended quadratic algebra and a model of the equivariant
cohomology ring of flag varieties}
\date{}
\author{Anatol N. Kirillov and Toshiaki Maeno}
\begin{document}
\maketitle
\begin{abstract}
For the root system of type $A$ we introduce and studied a certain extension of 
the quadratic algebra invented by S. Fomin and the first author, to construct 
a model for the equivariant cohomology ring of the corresponding flag variety.
As an application of our construction we describe a generalization  of the 
equivariant Pieri rule for double Schubert polynomials. For a general finite 
Coxeter system we construct an extension of the corresponding 
Nichols-Woronowicz algebra. In the case of finite  crystallographic Coxeter 
systems we present a construction of extended Nichols-Woronowicz algebra
model for the equivariant cohomology of the corresponding flag variety.
\end{abstract}
\section{Introduction}
In the paper \cite{FK} S. Fomin and the first author have introduced and 
study a model for the cohomology ring of flag varieties of type $A$ as a 
commutative subalgebra generated by the so-called Dunkl 
elements in a certain noncommutative quadratic algebra $\E_n$. 
An advantage of the approach developed in \cite{FK} is that it admits
a simple generalization which is suitable for description of the quantun 
cohomology ring of flag varieties, as well as (quantum) Schubert polynomials.  
Constructions from the paper \cite{FK} have 
been generalized to other finite root systems by the authors in \cite{KM1}. 
One of the main constituents of the above constructions is the Dunkl elements.
The basic properties of the Dunkl elements are:

1) they are pairwise commuting; 

2) in the so-called Calogero-Moser representation \cite{FK,KM1} they 
correspond to the {\it truncated} (i.e. without differential 
part) rational Dunkl operators \cite{D}; 

3) in the crystallographic case they correspond -- 
after applying the so-called Bruhat representation \cite{FK,KM1}-- to the Monk 
formula in the cohomology ring of the flag variety in question; 

4) in the crystallographic case, subtraction-free expressions 
of Schubert polynomials calculated at the Dunkl elements in the algebra 
$\widetilde{\cal{BE}}(\Delta)$ should provide a combinatorial rule for 
describing the Schubert basis structural constants, i.e. the intersection 
numbers of Schubert classes. 

In the case of classical root systems $\Delta,$ the 
first author \cite{Ki} has defined a certain extension 
$\widetilde{\cal{BE}}(\Delta)$ of the algebra ${\cal{BE}}(\Delta)$ together 
with a pairwise commuting family of elements, called Dunkl elements, which 
after applying the 
Calogero-Moser representation exactly coincide with the rational Dunkl 
operators. One of the main
objective of our paper is to study a commutative subalgebra generated by
the Dunkl elements in the extended algebra $\widetilde{\cal{BE}}(\Delta)$ in 
the 
case of type $A$ root systems. 
Postnikov \cite{Po} proved a Pieri-type formula for the elementary symmetric polynomial 
evaluated at the Dunkl elements, which was originally conjectured in \cite{FK}. 
In Theorem 2.1 we show an analogue of the Pieri-type formula in the extended algebra 
$\widetilde{\cal{BE}}(\Delta)$ of type $A.$ As a consequence, it is shown that 
the commutative subalgebra generated by the Dunkl elements is isomorphic 
to the equivariant cohomology ring $H^*_T(Fl_n),$ where $Fl_n$ is the flag variety 
parametrizing the full flags in the vector space ${\bf C}^n$ and $T=({\bf C}^{\times})^n$ 
is the torus acting on ${\bf C}^n$ diagonally. 

In Section 3 we construct the Bruhat representation of the 
algebra $\E_n\langle R \rangle[t]$ and study some properties of the former. 
The existence of Bruhat's representation of the algebra 
$\E_n\langle R \rangle [t]$ plays a crucial role in applications to the 
equivariant Schubert calculus, and constitutes an important step in the construction 
of the model of $H_T^*(Fl_n).$ Our formula in Theorem 2.1 describes the Pieri formula 
for the equivariant cohomology ring $H^*_T(Fl_n)$ via the Bruhat representation. 
As shown in Section 4, our formula also covers the quantized case. 
The Pieri formula for double Schubert polynomials was studied in \cite[Chapter 5]{V}. 
A similar rule in the cohomology ring $H_T^*(Fl_n)$ was stated and proved 
in \cite{R}. We will show how to compute the multiplication by the special Schubert classes 
via the Bruhat representation in Example 3.1, where an example is also given to see 
that the result of the computation based on our formula coincides with 
the one from \cite{R}. 

Another objective of our paper is to construct a certain extension of the
Nichols-Woronowicz model for the coinvariant algebra of a finite Coxeter
group $W$. 
It is conjectured that the Fomin-Kirllov quadratic algebra $\E_n$ is isomorphic 
to the Nichols-Woronowicz algebra associated to a certain kind of Yetter-Drinfeld 
module defined by the data of the root system of type $A_{n-1}.$ 
More generally, the algebra ${\cal{BE}}(\Delta)$ is a lift of the Nichols-Woronowicz 
algebra for the root system $\Delta.$ 
Recall that the Nichols-Woronowicz algebra model for the 
cohomology ring of flag varieties has been invented by Y. Bazlov \cite{Ba}. 
In Section 4 we introduce a certain extension $\tb_W$ of the Nichols-Woronowicz 
algebra $\B_{W}$ and construct a commutative subalgebra in the extended 
Nichols-Woronowicz algebra. Our second main result states that, 
for crystallographic root systems and $t=0,$ the commutative 
subalgebra of $\tb_W$ in question is isomorphic to the $T$-equivariant cohomology 
ring of the corresponding flag variety.

\section{Extension of the quadratic algebra}
\begin{dfn} 
{\it 
The algebra $\E_n$ is an associative algebra generated by the
symbols $[i,j],$ $1\leq i,j \leq n,$ $i\not=j,$ subject to
the relations: \\
$(0):$ $[i,j]=-[j,i]$ \\
$(1):$ $[i,j]^2=0,$ \\
$(2):$ $[i,j][k,l]=[k,l][i,j],$ if $\{ i,j \} \cap \{ k,l \} =\emptyset,$ \\
$(3):$ $[i,j][j,k]+[j,k][k,i]+[k,i][i,j]=0.$}
\end{dfn}
Let us consider the extension $\E_n\langle R \rangle [t]$ of the quadratic 
algebra
$\E_n$ by the
polynomial ring $R[t]=\Z[x_1,\ldots,x_n][t]$ defined by the commutation 
relations:
\\
(A): $[i,j]x_k= x_k[i,j],$ for $k\not=i,j,$ \\
(B): $[i,j]x_i=x_j[i,j]+t,$ $[i,j]x_j=x_i[i,j]-t,$ for $i<j,$ \\
(C): $[i,j]t=t[i,j].$

Note that the ${\bf S}_n$-invariant subalgebra $R^{{\bf S}_n}[t]$ of $R[t]$ 
is contained in the center of the algebra
$\E_n\langle R \rangle [t].$ Let $e_k(x_1,\ldots,x_n)$,~$1 \le k \le n,$ 
stands for the elementary symmetric 
polynomial of degree $k$ in the variables $x_1,\ldots,x_n.$ We put by 
definition, $e_{0}(x_1,\ldots,x_n)=1,$ and $e_k(x_1,\ldots,x_n)=0,$ if
$k < 0.$ 
\begin{dfn} {\it 
$(1)$ We define the $R[t]$-algebra $\te_n[t]$ by
\[ \te_n[t]= \E_n\langle R \rangle [t] \otimes_{R^{{\bf S}_n}} R . \]
More explicitly, $\te_n[t]$ is an algebra over the polynomial ring
$\Z[y_1,\ldots,y_n]$ generated by the symbols
$[i,j],$ $1\leq i,j \leq n,$ $i\not=j,$ and $x_1,\ldots,x_n,t$ satisfying
the relations in the definition of the algebra $\E_n\langle R \rangle [t],$
together with the identification $e_i(x_1,\ldots,x_n)=e_i(y_1,\ldots,y_n),$
for $i=1,\ldots,n.$ 
Denote by $\te_{n,t_0}$ the specialization of $\te_n [t]$ at $t=t_0.$
\\
$(2)$ The Dunkl elements $\theta_i\in \te_n [t],$ $i=1,\ldots,n,$ are
defined by the formula
\[ \theta_i=\theta_i^{(n)} = x_i + \sum_{j\not=i}[i,j] . \] }
\end{dfn}
\begin{rem}
{\rm Note that $x_i$'s do not commute with the Dunkl elements, but only
symmetric polynomials in $x_i$'s do. 
Moreover, we need the $R$-algebra structure of $\te_n[t]$ 
to construct the model of the 
$T$-equivariant cohomology ring $H^*_T(Fl_n)$ which is an algebra over 
$H^*_T({\rm pt}.) \cong R.$  
By this reason we need the second copy of
$R=\Z[y_1,\ldots,y_n]$, where $y_i$'s are assumed to belong to the center of 
the algebra $\te_n [t],$ and $f(x_1,\ldots,x_n)=f(y_1,\ldots,y_n)$ for any 
symmetric polynomial $f.$ }
\end{rem}
\begin{lem} {\it 
The Dunkl elements commute each other.}
\end{lem}
Proof.\quad This follows from the fact that
\[ (x_i+x_j)[i,j]=[i,j](x_i+x_j) . \]

\begin{thm} {\rm (Pieri formula in the algebra $\E_n\langle R \rangle [t]$)}
{\it For $k\leq m \leq n,$ we have 
$$e_k(\theta_1^{(n)},\ldots,\theta_m^{(n)})=$$
$$ \sum_{r \ge 0} (-t)^{r}~N(m-k,2r) \bigg\{ \sum_{S,I=\{ i_a \},(j_a)}
X_{S} \cdot 
[i_1,j_1] \cdots [i_{|I|},j_{|I|} ]\bigg\} ,$$
where 
\[ N(a,2b)= (2b-1)!! {a+2b \choose 2b} , \]
$X_{S}:= \prod_{s \in S}x_{s},$ and the second summation runs over triples
$(S,I=\{ i_1,\ldots,i_{|I|} \},(j_a)_{a=1}^{|I|})$ such that 
$S \subset \{ 1,\ldots, m \};$ $I$ is a subset of 
$\{1,\ldots, m \} \setminus S;$ 
$|I|+|S|+2r=k;$ $1 \le i_a \le m < j_a \le n;$ 
$j_1\leq \cdots \leq j_{|I|}.$} 
\end{thm}
Proof.\quad Let $\A$ be a subset of $\{ 1,\ldots,n \},$ $m:=|\A|,$ $d:=n-m$ and 
$\{ 1,\ldots,n \}\setminus \A =\{j_1<\cdots < j_d\}.$ 
Denote by $E_k(\A)$ the right-hand side of the formula, i.e., 
\[ E_k(\A):= \sum_{r\geq 0}(-t)^r N(m-k,2r)\sum_{S\subset \A}X_S\sum_{(*)} [s_1,t_1] 
\cdots [s_{k-2r-|S|},t_{k-2r-|S|}], \] 
where $(*)$ stands for the conditions that $s_1,\ldots, s_{k-2r-|S|}\in \A \setminus S$ 
are distinct, $t_1,\ldots,t_{k-2r-|S|} \in \{ 1,\ldots ,n \} \setminus \A$ and 
$t_1\leq \cdots \leq t_{k-2r-|S|}.$ 
It will 
suffice to prove the recursive formula 
\[ E_k(\A \cup \{ j=j_1 \})=E_k(\A)+E_{k-1}(\A)(x_j+\sum_{s\not= j}[j,s]). \]
For a subset $I=\{i_1,\ldots,i_l \} \subset \{ 1,\ldots, n \}$ and 
$p\not\in I,$ we use the symbol 
\[ \lla I | p \rra = \sum_{w\in \SS_l} [i_{w(1)},p]\cdots [i_{w(l)},p] \] 
as defined in \cite{Po}. 
We also use the symbol 
$I_1\cdots I_d \subset_m I$ which means 
that $I_1,\ldots,I_d \subset I$ are disjoint and 
$\# I_1 + \cdots + \# I_d=m.$ 
We have the following decompositions: 
\begin{eqnarray*} 
\lefteqn{E_k(\A)} \\ 
& = & \sum_{r \geq 0} (-t)^r N(m-k,2r) \sum_{S\subset \A}X_S 
\sum_{I_1\cdots I_d \subset_{k-2r-|S|} \A \setminus S} \!\!\!\!\!\!\! 
\lla I_1 |j_1 \rra \cdots \lla I_d | j_d \rra \\ 
& = & \!\! \sum_{r \geq 0}(-t)^r N(m-k,2r) (A_1^r+A_2^r), \\ 
\lefteqn{E_k(\A \cup \{ j \})} \\ 
& = & \sum_{r \geq 0} (-t)^r N(m-k+1,2r) \sum_{S\subset \A} X_S 
\sum_{I_2\cdots I_d \subset_{k-2r-|S|} \A \cup \{ j \} \setminus S} \!\!\!\!\!\!\!\!\!\!\!\!
\lla I_2 |j_2 \rra \cdots \lla I_d | j_d \rra \\ 
& = & \!\! \sum_{r \geq 0}(-t)^r N(m-k+1,2r) (B_1^r+B_2^r+B_3^r), \\ 
\lefteqn{E_{k-1}(\A)\sum_{s\not= j}[j,s]} \\ 
& = & \sum_{r\geq 0} (-t)^r N(m-k+1,2r) \sum_{S\subset \A}X_S 
\sum_{I_1\cdots I_d \subset_{k-1-2r-|S|} \A \setminus S} \!\!\!\!\!\!\! 
\lla I_1 |j_1 \rra \cdots \lla I_d | j_d \rra \sum_{s\not= j}[j,s] \\ 
& = & \sum_{r\geq 0} (-t)^r N(m-k+1,2r)  
(C_1^r+C_2^r+C_3^r+C_4^r), \\
\lefteqn{E_{k-1}(\A)x_j} \\ 
& = & \sum_{r\geq 0} (-t)^r N(m-k+1,2r) \sum_{S\subset \A}X_S 
\sum_{I_1\cdots I_d \subset_{k-1-2r-|S|} \A \setminus S} \!\!\!\!\!\!\! 
\lla I_1 |j_1 \rra \cdots \lla I_d | j_d \rra x_j \\ 
& =& \sum_{r\geq 0} (-t)^r N(m-k+1,2r)  
(D_1^r+D_2^r), 
\end{eqnarray*}
where $A_i^r,$ $B_i^r,$ $C_i^r,$ $D_i^r$ are defined as follows. 
\begin{itemize}
\item $A_1^r$ is the sum of terms with $I_1=\emptyset$; $A_2^r$ is the sum of
terms with $I_1\not= \emptyset$. 
\item $B_1^r$ is the sum of terms with $j\not\in S \cup I_2 \cup \cdots 
\cup I_d$; $B_2^r$ is the sum of terms with $j\in I_2 \cup \cdots \cup I_d$; 
$B_3^r$ is the sum of terms with $j\in S$. 
\item $C_1^r$ is the sum of terms with $s\in \A \setminus (S\cup I_1 \cup \cdots \cup I_d);$ 
$C_2^r$ is the sum of terms with $s\in I_2 \cup \cdots \cup I_d \cup \A^c$; 
$C_3^r$ is the sum of terms with $s\in S$; $C_4^r$ is the sum of 
terms with $s\in I_1.$ 
\item $D_1^r$ is the sum of terms with $I_1=\emptyset;$ $D_2^r$ is the sum of 
terms with $I_1\not= \emptyset.$  
\end{itemize}
Based on the same arguments used in \cite{Po}, we can see that 
$A_1^r=B_1^r,$ $A_2^r+C_1^r=0,$ ~$B_2^r=C_2^r$ and $C_4^r=0.$  It is also 
easy to see that $B_3^r=D_1^r.$  Now we have 
\begin{eqnarray*} 
\lefteqn{E_k(\A)+E_{k-1}(\A)(x_j+\sum_{s\not= j} [j,s])-E_k(\A \cup \{ j \})} \\ 
& = & \sum_{r\geq 0} (-t)^r \left( N(m-k,2r)(A_1^r+A_2^r)-
N(m-k+1,2r)(B_1^r-C_1^r-C_3^r-D_2^r) \right) \\ 
& = & \sum_{r\geq 1} (-t)^r \left( N(m-k,2r) - N(m-k+1,2r) \right) 
(A_1^r+A_2^r) \\ 
& & + \sum_{r\geq 0} (-t)^r N(m-k+1,2r)(C_3^r+D_2^r) . 
\end{eqnarray*} From the commutation relation $[i,j]x_j=x_i[i,j]-t,$ we have 
\begin{eqnarray*}
D_2^r  
& = & \sum_{S\subset \A}X_S \!\!\! \sum_{\scriptstyle I_1\cdots I_d \subset_{k-1-2r-|S|} \A \setminus S \atop 
I_1=\{a_1,\ldots , a_{|I_1|} \} } \sum_{w\in \SS_{|I_1|}} \!\! 
x_{a_{w(|I_1|)}} [a_{w(1)}, j] \cdots [a_{w(|I_1|)}, j] 
\lla I_2 |j_2 \rra \cdots \lla I_d | j_d \rra \\ 
& & -t \sum_{S\subset \A}X_S \!\!\!\sum_{I_1\cdots I_d \subset_{k-1-2r-|S|} \A \setminus S \atop 
I_1=\{a_1,\ldots , a_{|I_1|} \} } \sum_{w\in \SS_{|I_1|}} \!\! 
[a_{w(1)}, j] \cdots [a_{w(|I_1|-1)}, j] 
\lla I_2 |j_2 \rra \cdots \lla I_d | j_d \rra \\ 
&=& \sum_{S\subset \A}\sum_{s\not\in S}X_{S\cup \{ s \}} \sum_{I_1\cdots I_d \subset_{k-1-2r-(|S|+1)} \A \setminus S\cup \{ s \}} \lla I_1|j_1 \rra [s,j] 
\lla I_2 |j_2 \rra \cdots \lla I_d | j_d \rra \\ 
& & -(m-k+2r+2)t \sum_{S\subset \A} X_S 
\sum_{I_1\cdots I_d \subset_{k-2-2r-|S|} \A \setminus S}
\lla I_1 |j_1 \rra \cdots \lla I_d | j_d \rra \\ 
& = & -C_3^r +(-t)(m-k+2r+2)(A_1^{r+1}+A_2^{r+1}). 
\end{eqnarray*}
Hence, we have 
\begin{eqnarray*} 
\lefteqn{(-t)^{r+1} \left(N(m-k,2(r+1))-N(m-k+1,2(r+1) \right)(A_1^{r+1}+A_2^{r+1})} \\ 
& = & -(-t)^{r+1} (2r+1)!! \frac{(m-k+2r+2)!}{(2r+1)!(m-k+1)!}(A_1^{r+1}+A_2^{r+1}) \\ 
& = & -(-t)^r(2r-1)!! \frac{(m-k+2r+1)!}{(2r)!(m-k+1)!}\cdot 
(-t)(m-k+2r+2)(A_1^{r+1}+A_2^{r+1}) \\ 
& = & -(-t)^r N(m-k+1,2r)(C_3^r+D_2^r). 
\end{eqnarray*}
This shows the desired result. 

\begin{ex}
{\rm Let us check the coefficients of $t$ and $t^2$ in the expression of 
$\theta_1\theta_2\theta_3\theta_4 \in \E_5\langle R \rangle [t]$ by direct computation. 
It is easy to see 
\begin{eqnarray*}
\theta_1 \theta_2 \theta_3 
& = & -t(x_1+x_2+x_3+[14]+[15]+[24]+[25]+[34]+[35]) \\ 
& & +x_1x_2x_3+x_1x_2([34]+[35])+x_1x_3([24]+[25])+x_2x_3([14]+[15]) \\ 
& & + x_1 ([24][34]+[24][35]+[25][35]+[34][24]+[34][25]+[35][25]) \\ 
& & + x_2 ([14][34]+[14][35]+[15][35]+[34][14]+[34][15]+[35][15]) \\ 
& & + x_3 ([14][24]+[14][25]+[15][25]+[24][14]+[24][15]+[25][15]) \\ 
& & +\sum_{ \{i_1,i_2,i_3\}=\{1,2,3 \},\; 4\leq j_1\leq j_2 \leq j_3\leq 5}[i_1j_1][i_2j_2][i_3j_3]. 
\end{eqnarray*} 
Multiply this expression by $\theta_4$ from the right. We have the contributions 
to the coefficients of $t$ and $t^2$ only from the following terms: 
\begin{eqnarray*} 
\lefteqn{\bullet \;\;\;\; -t \Big( x_1+x_2+x_3+[14]+[15]+[24]+[25]+[34]+[35] \Big) \theta_4} \\ 
& = & 3t^2-t \Big( x_1x_4 + x_2x_4+x_3x_4+ x_1([42]+[43]+[45]) \\  
& & +x_2([41]+[43]+[45])+x_3([41]+[42]+[45])+x_4([15]+[25]+[35]) \\ 
& & + [15][41]+[24][41]+[25][41]+[34][41]+[35][41] \\ 
& & + [14][42]+[15][42]+[25][42]+[34][42]+[35][42] \\ 
& & + [14][43]+[15][43]+[24][43]+[25][43]+[35][43] \\ 
& & + [14][45]+[15][45]+[24][45]+[25][45]+[34][45]+[35][45] \Big), \\
 \\
\lefteqn{\bullet \;\;\;\; x_1 ([24][34]+[24][35]+[25][35]+[34][24]+[34][25]+[35][25])x_4} \\ 
&=& -tx_1([24]+[34]+[25]+[35]) + \cdots , \\ 
\lefteqn{\bullet \;\;\;\; x_2 ([14][34]+[14][35]+[15][35]+[34][14]+[34][15]+[35][15])x_4} \\ 
&=& -tx_2([14]+[34]+[15]+[35]) + \cdots , \\ 
\lefteqn{\bullet \;\;\;\; x_3 ([24][14]+[24][15]+[25][15]+[14][24]+[14][25]+[15][25])x_4} \\ 
&=& -tx_3([24]+[14]+[25]+[15]) +\cdots ,  
\end{eqnarray*}
\begin{eqnarray*}
\bullet \;\;\;\; x_1x_2([34]+[35])x_4 &=& -tx_1x_2+x_1x_2(x_3[34]+x_4[35]), \;\;\;\;\;\;\;\;\;\;\;\;\;\;\;\;\;\;\;\; \\ 
\bullet \;\;\;\; x_2x_3([14]+[15])x_4 &=& -tx_2x_3+x_2x_3(x_1[14]+x_4[15]), \;\;\;\;\;\;\;\;\;\;\;\;\;\;\;\;\;\;\;\; \\ 
\bullet \;\;\;\; x_1x_3([24]+[25])x_4 &=& -tx_1x_3+x_1x_3(x_2[24]+x_4[25]), \;\;\;\;\;\;\;\;\;\;\;\;\;\;\;\;\;\;\;\;
\end{eqnarray*}
\begin{eqnarray*}
\lefteqn{\bullet \;\;\;\;\;\;\;\;\;\; \sum_{ \{i_1,i_2,i_3\}=\{1,2,3 \},\; 4\leq j_1\leq j_2 \leq j_3\leq 5}
[i_1j_1][i_2j_2][i_3j_3] x_4} \\ 
&=& -t([14][24]+[14][35]+[25][35]+[24][14]+[24][35]+[15][35] \\ 
& & +[34][14]+[34][25]+[15][25]+[34][24]+[34][15]+[25][15] \\ 
& & +[14][34]+[14][25]+[35][25]+[24][34]+[24][15]+[35][15])+\cdots . 
\end{eqnarray*}
By using the relations $[15][41]+[14][45]=[45][15],$ $[25][42]+[24][45]=[45][25],$ 
$[35][43]+[34][45]=[45][35],$ we obtain finally that 
\begin{eqnarray*}
\theta_1\theta_2\theta_3\theta_4 
&=& 3t^2-t \Big( x_1x_2+x_1x_3+x_1x_4+x_2x_3+x_2x_4+x_3x_4 \\ 
& & + x_1([25]+[35]+[45])+x_2([15]+[35]+[45]) \\ 
& & + x_3([15]+[25]+[45])+x_4([15]+[25]+[35]) \\ 
& & +[15][25]+[15][35]+[15][45]+[25][35]+[25][45]+[35][45] \\ 
& & +[25][15]+[35][15]+[45][15]+[35][25]+[45][25]+[45][35]\Big) + \cdots . 
\end{eqnarray*} 
It is easy to see that the formula for $\theta_1\theta_2\theta_3\theta_4$ stated 
in Theorem 2.1 produces the same expression. }
\end{ex}

The following is a special case of the formula in Theorem 2.1 for $m=n.$ 
\begin{cor} {\it We have the relations 
$$e_k(\theta_1^{(n)},\ldots,\theta_n^{(n)})=$$
$$ e_k(y_1,\ldots,y_n)+
\sum_{r \ge 1}(-t)^{r}~
(2r-1) !!~{n-k+2r \choose 2r}~e_{k-2r}(y_1,\ldots,y_n), \; 1 \le k \le n, $$
in the algebra $\te_n [t].$} 
\end{cor}
It will be shown in Section 3 that the above relations describe the complete set of 
relations among the Dunkl elements in $\te_n [t].$

\section{Bruhat representation}
Let us recall the definition of the Bruhat representation
of the algebra $\E_n$ on the group ring of the symmetric group
$\Z \langle {\bf S}_n \rangle =
\oplus_{w\in {\bf S}_n} \Z  \cdot \underline{w}.$
The operator $\sigma_{ij},$ $i<j,$ is defined as follows:
\[ \sigma_{ij} ( \underline{w}) = \left\{ \begin{array}{cc}
\underline{wt_{ij}}, & \textrm{if $l(wt_{ij})=l(w)+1,$} \\
0, & \textrm{otherwise,}
\end{array} \right. \]
where $t_{ij} \in {\bf S}_n$ is the transposition of $i$ and $j.$ 
Then the Bruhat representation of $\E_n$ is defined by $[i,j].\underline{w}:=
\sigma_{ij}(\underline{w}).$

Now we extend the Bruhat representation to that of the algebra
$\E_n\langle R \rangle [t]$ defined on
\[ R[t] \langle {\bf S}_n \rangle = \oplus_{w\in {\bf S}_n} \Z[y_1,\ldots ,y_n][t] \cdot
\underline{w}. \]
Let us define the devided difference operator $\partial_{ij}$ on $\Z[y_1,\ldots ,y_n][t]$ 
as a $\Z[t]$-linear operator given by 
$\partial_{ij}:= (1-t_{ij})/(y_i-y_j).$ 
For $f(y) \in \Z[y_1,\ldots ,y_n][t]$ and $w\in {\bf S}_n,$ we define the
$\Z[t]$-linear operators
$\tilde{\sigma}_{ij},$ $i<j,$ and $\xi_k$ as follows:
\[ \tilde{\sigma}_{ij} (f(y)\underline{w})= \left\{ \begin{array}{cc}
t(\pa_{w(i)w(j)} f(y))\underline{w}+
f(y)\underline{wt_{ij}}, & \textrm{if $l(wt_{ij})=l(w)+1,$} \\
t(\pa_{w(i)w(j)} f(y))\underline{w}, & \textrm{otherwise,}
\end{array} \right. \]
\[ \xi_k  (f(y)\underline{w}) = (y_{w(k)}f(y))\underline{w} . \]
\begin{prop} {\it 
The algebra $\E_n\langle R \rangle [t]$ acts $\Z[t]$-linearly on $\Z[y][t] \langle {\bf S}_n \rangle$ via
$[ij] \mapsto \tilde{\sigma}_{ij}$ and $x_k \mapsto \xi_k.$}
\end{prop}
Proof.\quad Let us check the compatibility with the defining relations
of the algebra $\te_n [t].$ We show the compatibility only with the relations
(1), (3) and (B). The rest are easy to check.

Let us start with the relation (1). We have
\begin{eqnarray*} 
\tilde{\sigma}_{ij}^2 (f(y)\underline{w}) & = & 
\tilde{\sigma}_{ij} \left( t(\pa_{w(i)w(j)} f(y))\underline{w}+
f(y)\sigma_{ij}(\underline{w}) \right) \\
& = & t^2 (\pa_{w(i)w(j)}^2 f(y))\underline{w} +
t (\pa_{w(i)w(j)} f(y)) \sigma_{ij}(\underline{w}) \\
& & + t(\pa_{w(j)w(i)} f(y)) \sigma_{ij}(\underline{w}) +
f(y)\sigma_{ij}^2(\underline{w}) . 
\end{eqnarray*}
Since $\pa_{w(i)w(j)}^2=0,$ $\sigma_{ij}^2=0$ and $\pa_{w(i)w(j)}=-\pa_{w(j)w(i)},$
we get $\tilde{\sigma}_{ij}^2=0.$

For the relation (3), we have
\begin{eqnarray*} 
\tilde{\sigma}_{ij}\tilde{\sigma}_{jk} (f(y)\underline{w}) & = & 
\tilde{\sigma}_{ij} \left( t(\pa_{w(j)w(k)} f(y))\underline{w}+
f(y)\sigma_{jk}(\underline{w}) \right) \\
& = & t^2 (\pa_{w(i)w(j)} \pa_{w(j)w(k)} f(y))\underline{w} +
t (\pa_{w(j)w(k)} f(y)) \sigma_{ij}(\underline{w}) \\ 
& & + t(\pa_{w(i)w(k)} f(y)) \sigma_{jk}(\underline{w}) +
f(y)\sigma_{ij}\sigma_{jk}(\underline{w}) . 
\end{eqnarray*}
We also obtain $\tilde{\sigma}_{jk}\tilde{\sigma}_{ki} (f(y)\underline{w})$
and $\tilde{\sigma}_{ki}\tilde{\sigma}_{ij} (f(y)\underline{w})$ by the cyclic
permutation of $i,j,k.$ The 3-term relations
\[ \pa_{w(i)w(j)} \pa_{w(j)w(k)}+\pa_{w(j)w(k)} \pa_{w(k)w(i)}+
\pa_{w(k)w(i)} \pa_{w(i)w(j)} =0 \]
and
\[ \sigma_{ij}\sigma_{jk} + \sigma_{jk}\sigma_{ki} + \sigma_{ki}\sigma_{ij}=0 \]
show the desired equality
\[ \tilde{\sigma}_{ij}\tilde{\sigma}_{jk} + \tilde{\sigma}_{jk}\tilde{\sigma}_{ki} +
\tilde{\sigma}_{ki}\tilde{\sigma}_{ij} =0. \]
Finally, we check the relation (B). We have
\begin{eqnarray*} 
\tilde{\sigma}_{ij}\xi_i(f(y)\underline{w}) & = & 
\tilde{\sigma}_{ij}(y_{w(i)}f(y)\underline{w}) \\
& = & t\pa_{w(i)w(j)}(y_{w(i)}f(y))\underline{w}+
(y_{w(i)}f(y))\sigma_{ij}(\underline{w}) \\
& = & t(f(y)\underline{w})+t(y_{w(j)}\pa_{w(i)w(j)}f(y))\underline{w} +
y_{wt_{ij}(j)}\sigma_{ij}(\underline{w}) \\ 
& = & \xi_j \tilde{\sigma}_{ij}(f(y)\underline{w}) + t (f(y)\underline{w}) . 
\end{eqnarray*}
Now the proof of the well-definedness of the Bruhat representation is completed. 

Let us denote by $AH_n^{0}$ the subalgebra in 
$\E_n\langle R \rangle [t]$ generated by the elements $t,$ 
$h_1:=[1,2],h_2:=[2,3], \ldots,$ $h_{n-1}:=[n-1,n]$ and 
$x_1,\ldots,x_n.$ 
\begin{lem} {\it 
The module $R[t] \langle {\bf S}_n \rangle$ is generated by $\underline{\rm{id.}}$ 
over $AH_n^0.$} 
\end{lem}
Proof.\quad If we take a reduced decomposition $w=s_{i_1}\cdots s_{i_l},$ $s_i=(i,i+1),$ 
of a permutation $w\in {\bf S}_n,$ then we have 
\[ h_{i_l}\cdots h_{i_1}(\underline{\rm{id.}})=\underline{w}. \] 
This shows $R[t] \langle {\bf S}_n \rangle = AH_n^0 \cdot \underline{\rm{id.}}$ 

\begin{thm} {\it The subalgebra $AH_n^0$ of $\E_n\langle R \rangle [t]$ 
is isomorphic to the nil degenerate affine Hecke algebra ${\cal AH}_n^0$ 
of type $A_{n-1}^{(1)},$ i.e. the $\Z[t]$-algebra given by two sets of 
generators $g_1,\ldots,g_{n-1}$ and ${x_1,\ldots,x_n}$ subject to the set of 
defining relations:} {\rm 
$$ g_i^2=0, ~~g_ig_j=g_jg_i, \; \; \textrm{if $|i-j| > 1$},~~g_ig_jg_i=
g_jg_ig_j,~\textrm{if $|i-j|=1$},$$
$$ x_ix_j=x_jx_i, ~~x_kg_i=g_ix_k,~\textrm{if $k \not= i,i+1,$}~~
g_ix_i-x_{i+1}g_i=t.$$}
\end{thm}
Proof.\quad First of all it is easy to see that the elements 
$t, h_1,\ldots,h_{n-1}, x_1,\ldots,x_n$ do satisfy the relations listed above. 
Hence we have a surjective homomorphism $\rho$ from the nil degenerate affine Hecke 
algebra ${\cal AH}_n^0$ to $AH_n^0$ given by $\rho(g_i)=h_i.$ 
Now we are going to 
construct a basis in the algebra ${\cal AH}_n^{0}.$ Let $w \in {\bf S}_n$ be a 
permutation and $w=s_{i_1} \cdots s_{i_k}$ its any reduced decomposition. 
Since the elements $g_1,\ldots,g_{n-1}$ satisfy the Coxeter relations, the element 
$g_{w}:=g_{i_1} \cdots g_{i_k}$ is well-defined. On the other hand, using 
relations among the elements $x_1, \ldots,x_n$ and $g_1, \ldots,g_{n-1},$ 
in the algebra ${\cal AH}_n^{0},$ one can write any element of 
${\cal AH}_n^{0}$ as a linear combination of elements $t^ax^{m}g_{w},$ 
where $w \in {\bf S}_n,$ 
$m \in \Z_{\ge 0}^{n}$ and $a \in \Z_{\ge 0}.$ 
Lemma 3.1 implies that $R[t] \langle {\bf S}_n \rangle = 
\rho({\cal AH}_n^0) \cdot \underline{\rm{id.}},$ and hence the elements 
$\{t^ax^mg_w~|~a \in \Z_{\geq 0},~w \in {\bf S}_n,~m \in \Z_{\ge 0}^{n} \}$ 
must be linearly independent over ${\bf Z}$. 
This means the injectivity of $\rho,$ so we conclude that the homomorphism $\rho$ is 
isomorphism. \bigskip 

Let us consider the double Schubert polynomials $\S_w(x,y),$ $w\in {\bf S}_n,$ 
introduced by Lascoux and Sch\"utzenberger \cite{LS}. The polynomial $\S_{w_0}(x,y)$ 
for the maximal element $w_0\in {\bf S}_n$ is by definition given by 
$\S_{w_0}(x,y):=\prod_{i+j\leq n}(x_i-y_j).$ For $w\in {\bf S}_n,$ take a reduced 
decomposition $w=s_{i_1}\cdots s_{i_l}.$ The divided difference operator 
$\partial_w:=\partial_{i_1}\cdots \partial_{i_l}$ is well-defined thanks to the 
Coxeter relations. We define the double Schubert polynomial for $w$ by $\S_w(x,y):=
\partial_{w^{-1}w_0}^{(x)}\S_{w_0}(x,y),$ where $\partial_w^{(x)}$ means a divided 
difference operator on $x$-variables. The Schubert polynomials 
are defined to be the specializations $\S_w(x):=\S_w(x,0)$ of the double 
Schubert polynomials. 

\begin{thm} {\it 
Let $\S_w(x,y)$ be the double Schubert polynomial corresponding to $w\in {\bf S}_n.$
When $t=0$, we have
\[ \S_w(\theta,y)(\underline{\rm{id.}})=\underline{w}. \]}
\end{thm}
Proof.\quad This follows from the Monk formula for the double Schubert
polynomials, see e.g. \cite[Exercise 2.7.2]{M}, and
\begin{eqnarray*} 
(\theta_i-y_{w(i)})(\underline{w}) & = &
\xi_i(\underline{w})+\sum_{j\not=i}\sigma_{ij}(\underline{w})-
y_{w(i)}\underline{w} \\ 
& = & \sum_{j>i,l(wt_{ij})=l(w)+1}\underline{wt_{ij}}-
\sum_{j<i,l(wt_{ij})=l(w)+1}\underline{wt_{ij}}. 
\end{eqnarray*}
Let $w\in {\bf S}_n,$ $r$ be the maximal descent of $w,$ and 
$s$ be the greatest integer such that $w(s)<w(r).$  
By applying the Monk formula to the product $x_r \S_u(x,y),$ $u=wt_{rs},$ 
we get the transition formula 
\[ \S_w(x,y)=(x_r-y_{u(r)})\S_u(x,y)+\sum_{v\in S(w,r)}\S_v(x,y), \] 
where $S(w,r)$ is the set of permutations of form $ut_{jr}=wt_{rs}t_{jr}$ 
with $j<r$ and of the same length as $w.$ 
Similarly we have 
\[ \underline{w}=(\theta_r-y_{u(r)})\underline{u}+\sum_{v\in S(w,r)}\underline{v}. \] 
Let $r'$ be the maximal descent of $v=wt_{rs}t_{jr}\in S(w,r).$ 
Since $l(v)=l(wt_{rs})+1,$ we have $v(r)=wt_{rs}(j)<wt_{rs}(r)=w(s)<w(r).$ 
Hence we have $r'<r$ or ``$r'=r$ and $v(r)<w(r)$''. 
Then our assertion is proved by induction on the Bruhat ordering, 
the maximal descent $r$ and $w(r).$ 

\begin{rem}
{\rm Only when $t=0,$ one can extend $\Z[y][t]$-linearly 
the Bruhat representation of the algebra 
$\E_n\langle R \rangle [t]$ to that of the algebra $\te_{n,0}.$ 
In this case, the Dunkl elements commutes with the multiplication 
by $y_i$'s. }
\end{rem}

\begin{prop} {\it 
The list of relations given in Corollary 2.1 
describes the complete set of relations among the Dunkl elements 
$\theta_1^{(n)},\ldots,\theta_n^{(n)}$ in the algebra $\te_n[t].$ 
In other words, the following surjective homomorphism $\varphi$ 
between $\Z[t][y_1,\ldots,y_n]$-algebras is an isomorphism: 
\[ \begin{array}{cccc}
\varphi : & \Z[t][y_1,\ldots,y_n][z_1,\ldots,z_n]/J_n^t & \rightarrow & 
\Z[t][y_1,\ldots,y_n][\theta_1,\ldots,\theta_n] \subset \te_n[t] \\ 
 & z_i & \mapsto & \theta_i, 
\end{array} \]
where the ideal $J_t$ is generated by the polynomials 
\[ e_k(z_1,\ldots,z_n) - e_k(y_1,\ldots,y_n) - 
\sum_{r \ge 1}(-t)^{r}~
(2r-1) !!~{n-k+2r \choose 2r}~e_{k-2r}(y_1,\ldots,y_n) \] 
for $k=1,\ldots,n.$ } 
\end{prop}
Proof.\quad 
For the (single) Schubert polynomial 
$\S_w(\theta_1,\ldots,\theta_n)\in \E_n \langle R \rangle [t]$ in the Dunkl elements, 
we have 
\[ \S_w(\theta_1,\ldots,\theta_n)(\underline{\rm{id}.})|_{t=0}=
\underline{w}+(\textrm{linear combination of $\underline{v}$ with $l(v)<l(w)$}), \] 
so 
\begin{eqnarray*} 
\lefteqn{\S_w(\theta_1,\ldots,\theta_n)(\underline{\rm{id}.})} \\ 
 & = & 
\S_w(\theta_1,\ldots,\theta_n)(\underline{\rm{id}.})|_{t=0}
+ (\textrm{linear combination of $\underline{v}$ with $l(v)<l(w)$}) \\ 
 & = & \underline{w}+ (\textrm{linear combination of $\underline{v}$ with $l(v)<l(w)$}). 
\end{eqnarray*} 
Hence, $\S_w(\theta_1,\ldots,\theta_n)$'s are linearly independent in 
$\E_n \langle R \rangle [t]$ over $R^{{\bf S}_n}.$ 

Let $R_{{\bf S}_n}$ be the coinvariant algebra of ${\bf S}_n.$ Since 
$R=R^{{\bf S}_n}\otimes R_{{\bf S}_n},$ the polynomials $\S_w(z),$ $w\in {\bf S}_n,$ 
form a $\Z[z]^{{\bf S}_n}$-basis of $\Z[z].$ In particular, any polynomial 
$f(z_1,\ldots,z_n) \in \Z[z]$ can be expressed as 
\[ f(z_1,\ldots,z_n)= \sum_{w\in {\bf S}_n}\phi_w(z)\S_w(z), \] 
where $\phi_w(z)\in \Z[z]^{{\bf S}_n}=\Z[e_1(z),\ldots,e_n(z)].$ 
Therefore the image of $f(z_1,\ldots,z_n)$ 
in $\Z[t][y][z_1,\ldots,z_n]/J_n^t$ is a linear combination of 
$\S_w(z)$'s over $\Z[t][y].$ Since the elements $\S_w(\theta_1,\ldots,\theta_n)$ 
in $\te_n[t]$ are linearly independent over $\Z[t][y],$ the homomorphism 
$\varphi$ is an isomorphism. 
\begin{cor} {\it
The subalgebra of $\te_{n,0}$ generated by the Dunkl elements
$\theta_1,\ldots,\theta_n$
over $H^*_T({\rm pt})=\Z[y_1,\ldots,y_n]$ is isomorphic
to the $T$-equivariant cohomology ring $H^*_T(Fl_n).$}
\end{cor}
Proof.\quad Let $(0=U_0 \subset U_1 \subset \cdots \subset U_n)$ be the universal 
flag over $Fl_n.$ 
First of all it follows from Corollary 2.1 that the natural map 
$z_i:=-c_1^T(U_i/U_{i-1}) \mapsto \theta_i,$ $y_i \mapsto y_i$ defines a
surjective homomorphism 
\[ \pi:H_{T}^{*}(Fl_n) \rightarrow 
{\Z}[y_1,\ldots,y_n][\theta_1,\ldots,\theta_n] \subset \te_{n,0} . \] 
On the other hand, Proposition 3.2 shows that the homomorphism $\pi$ is injective. 
\begin{ex} {\rm 
When $t=0,$ our formula in Theorem 2.1 specializes to 
the Pieri formula in $H^*_T(Fl_n)$ under the identification in Corollary 3.1. 
Let $\{ \Omega_w \}_{w\in {\bf S}_n}$ be the Schubert basis of $H^*_T(Fl_n)$ and 
$z_i:=-c_1^T(U_i/U_{i-1}).$ Then our formula in $H^*_T(Fl_n)$ can be written 
as follows: 
\[ e_k(z_1,\ldots,z_m)\cdot \Omega_w = \sum_{S\subset \{1,\ldots, m\} } \prod_{i\in S}y_{w(i)}
\sum_{(*)}\Omega_{wt_{i_lj_l}\cdots t_{i_1j_1}}, \]
where $(*)$ stands for the conditions $|S|+l=k;$ $1\leq i_a \leq m <j_a \leq n$ 
for $1\leq a \leq l;$ $i_1,\ldots,i_l$ are distinct; $j_1\leq \cdots \leq j_l;$ 
there exists a path $w \rightarrow wt_{i_lj_l} \rightarrow wt_{i_lj_l}t_{i_{l-1}j_{l-1}} 
\rightarrow \cdots \rightarrow wt_{i_lj_l}\cdots t_{i_1j_1}$ in the Bruhat ordering 
of ${\bf S}_n.$ 
For the cyclic permutation $[m,k]:=s_{m-k+1}s_{m-k+2} \cdots s_m,$ the corresponding 
double Schubert polynomial is given as follows (see \cite[Proposition 2.6.7]{M}): 
\[ \S_{[m,k]}(x,y)= \sum_{j=0}^ke_{k-j}(x_1,\ldots,x_m)h_j(-y_1,\ldots,-y_{m-k+1}), \] 
where $h_j$ is the complete symmetric polynomial of degree $j,$  
so the above formula 
gives the multiplication rule for $\Omega_{[m,k]}$ in $H^*_T(Fl_n).$ 

Now we consider an example from \cite{R}. 
Let $[5,5]=s_1s_2s_3s_4s_5$ be a permutation in ${\bf S}_9.$ 
For $w=s_3s_4s_5s_7s_6s_5$ and 
$u=wt_{26}t_{16}t_{59},$ let us compute the coefficient $p_{[5,5],w}^u$ of $\Omega_u$ in the 
expansion of $\Omega_{[5,5]}\cdot \Omega_w$ by using our formula. 
It is easy to see that there exists a unique path of length 3 from $w$ to $u$ as 
follows: 
\[ w \rightarrow wt_{59} \rightarrow wt_{59}t_{26}\rightarrow u=wt_{59}t_{26}t_{16}. \] 
Hence, the action of $\S_{[m,k]}(\theta,y)$ on $\underline{w}$ under the Bruhat 
representation at $t=0$ is given by  
\begin{eqnarray*}
\S_{[m,k]}(\theta,y)(\underline{w}) &=& \sum_{j=0}^5(-y_1)^je_{5-j}(\theta_1,\ldots,\theta_5) 
\underline{w} \\ 
 &=& (\cdots + (x_3x_4-y_1(x_3+x_4)+y_1^2)[16][26][59]+ \cdots )\underline{w} \\ 
 &=& \cdots + (y_1-y_4)(y_1-y_6)\underline{u}+ \cdots , 
\end{eqnarray*}
so we get $p_{[5,5],w}^u=(\alpha_1+\alpha_2+\alpha_3)(\alpha_1+\cdots + \alpha_5),$ 
$\alpha_i=y_i-y_{i+1}.$ This coincides with the result computed in \cite[Example 4.8]{R}.} 
\end{ex}
\section{Quantization}
\begin{dfn} {\it 
The algebra $\E_n^{\bf q}$ is a $\Z[q_{ij}=q_{ji}|1\leq i<j\leq n]$-algebra 
defined by the same generators and relations as in the definition 
of the algebra $\E_n$ except that the relation $(1)$ is 
replaced by 
\[ (1)' \;\;\; [i,j]^2=q_{ij} \] 
for $1\leq i<j \leq n.$ 
The algebra $\E_n^q$ is defined as a $\Z[q_1,\ldots,q_{n-1}]$-algebra 
obtained from $\E_n^{\bf q}$ by the specialization 
{\rm 
\[ q_{ij}= \left\{ 
\begin{array}{cc} 
q_i & \textrm{if $i=j-1,$} \\ 
0, & \textrm{if $i<j-1.$} 
\end{array}
\right. \] } 
The extension $\E_n^{\bf q} \langle R \rangle [t]$ (resp. $\E_n^q \langle R \rangle [t]$) 
of the algebra $\E_n^{\bf q}$ (resp. $\E_n^q$) 
is also defined by the relations $(A),$ $(B)$ and $(C).$} 
\end{dfn}
In the algebra $\E_n^{\bf q} \langle R \rangle [t],$ we have an analogous formula 
to Theorem 2.1. In order to state the formula, we need the quantum 
elementary symmetric polynomials $e_k^{\bf q}$ defined by the recursive formula 
\[ e_k^{\bf q}(X_i|i \in I\cup \{ j \})=e_k^{\bf q}(X_i|i \in I)+
X_je_{k-1}^{\bf q}(X_i|i \in I) + 
\sum_{a\in I}q_{aj}e_{k-2}^{\bf q}(X_i|i \in I\setminus \{ a \}), \] 
\[ e_0^{\bf q}(X_i|i \in I)=1, \; \; \; e_k^{\bf q}(\emptyset)=0, \; k> 0, \] 
where $I$ is a subset of $\{1,\ldots, n \}$ and $j\not\in I.$ 
The Dunkl elements $\theta_i=\theta^{(n)}_i$ in $\E_n^{\bf q} \langle R \rangle [t]$ 
are defined by the same formula in the classical case. 
\begin{thm} For $k\leq m \leq n,$ we have the following formula in the algebra 
$\E_n^{\bf q} \langle R \rangle [t]:$ 
$$e_k^{\bf q}(\theta_1^{(n)},\ldots,\theta_m^{(n)})=$$
$$ \sum_{r \ge 0} (-t)^{r}~N(m-k,2r) \bigg\{ \sum_{S,I=\{ i_a \},(j_a)}
X_{S} \cdot 
[i_1,j_1] \cdots [i_{|I|},j_{|I|} ]\bigg\} ,$$
where the second summation runs over triples
$(S,I=\{ i_1,\ldots,i_{|I|} \},(j_a)_{a=1}^{|I|})$ such that 
$S \subset \{ 1,\ldots, m \};$ $I$ is a subset of 
$\{1,\ldots, m \} \setminus S;$ 
$|I|+|S|+2r=k;$ $1 \le i_a \le m < j_a \le n;$ 
$j_1\leq \cdots \leq j_{|I|}.$
\end{thm} 
Proof.\quad We use the same symbols as in the proof of Theorem 2.1. 
We will show that $E_k(\A)$ satisfies the recursive relation 
\[ E_k(\A \cup \{ j \})=E_k(\A)+E_{k-1}(\A)(x_j+\sum_{s\not= j}[j,s]) + 
\sum_{\nu \in \A}q_{\nu j}E_{k-2}(\A\setminus \{ \nu \}) \] 
in the algebra $\E_n^{\bf q} \langle R \rangle [t].$ 
All the arguments in the proof of Theorem 2.1 work well except that $C_4^r=0,$ 
so 
\[ E_k(\A \cup \{ j \})-E_k(\A)-E_{k-1}(\A)(x_j+\sum_{s\not= j}[j,s]) = - 
\sum_{r\geq 0}(-t)^r N(m-k+1,2r) C_4^r. \] 
The following cyclic relation (\cite[Lemma 5.3]{Po}) holds in the algebra 
$\E_n^{\bf q} \langle R \rangle [t]:$ 
\begin{eqnarray*}
\lefteqn{\sum_{k=1}^m[a, i_k][a,i_{k+1}] \cdots [a, i_m]\cdot 
[a,i_1] \cdots [a,i_{k-1}][a,i_k]} \\ 
 & = & \sum_{k=1}^m 
q_{ak}[i_k,i_{k+1}][i_k,i_{k+2}]\cdots [i_k,i_m][i_k,i_1]\cdots [i_k,i_{k-1}], 
\end{eqnarray*}
where $1\leq a,i_1,\ldots,i_m \leq n$ are distinct. By using this cyclic 
relation, we have 
\begin{eqnarray*}
\lefteqn{-\sum_{r\geq 0}(-t)^r N(m-k+1,2r) C_4^r} \\ 
 & = & -\sum_{r\geq 0}(-t)^rN(m-k+1,2r) \sum_{S\subset \A}X_S \!\!
\sum_{I_1\cdots I_d \subset_{k-1-2r-|S|}\A\setminus S} \!\!\!\!\!\!\! \lla I_1|j_1 \rra \cdots 
\lla I_d |j_d \rra \sum_{s\in I_1}[j,s] \\ 
 & = & \sum_{r\geq 0}(-t)^rN(m-k+1,2r) \sum_{S\subset \A}X_S \sum_{\nu \in \A\setminus S}
\sum_{(*)} 
q_{j\nu}\lla I_1|j_1 \rra \cdots 
\lla I_d |j_d \rra  \lla I_{d+1} | \nu \rra \\ 
 & = & \sum_{\nu \in \A}q_{j\nu}\sum_{r\geq 0}(-t)^rN((m-1)-(k-2),2r) \!\!\!\! \sum_{S\subset \A\setminus \{ \nu \}}X_S 
\sum_{(*)} 
\lla I_1|j_1 \rra \cdots 
\lla I_{d+1} |j_{d+1} \rra  \\ 
 & = & \sum_{\nu \in \A}q_{j\nu}E_{k-2}(\A\setminus \{ \nu \}), 
\end{eqnarray*} 
where $(*)$ means the condition 
$I_1\cdots I_{d+1} \subset_{k-2-2r-|S|}\A \setminus (S\cup \{ \nu \}).$ 
This completes the proof. \medskip 

The Bruhat representation for $\E_n$ is deformed to the quantum 
Bruhat representation for $\E_n^q.$ 
We define the quantum Bruhat operator $\sigma_{ij}^q,$ $i<j,$ acting 
on $\Z[q_1,\ldots,q_{n-1}]\langle {\bf S}_n \rangle = 
\oplus_{w\in {\bf S}_n} \Z [q_1,\ldots,q_{n-1}] \cdot \underline{w}$ 
as follows: 
\[ \sigma_{ij}^q ( \underline{w}) = \left\{ \begin{array}{cc}
q_{ij} \underline{wt_{ij}}, & \textrm{if $l(wt_{ij})=l(w)-2(j-i)+1,$} \\
\underline{wt_{ij}}, & \textrm{if $l(wt_{ij})=l(w)+1,$} \\
0, & \textrm{otherwise.}
\end{array} \right. \]

For $f(y) \in \Z[y_1,\ldots ,y_n][t]$ and $w\in {\bf S}_n,$ we define the
$\Z[q_1,\ldots,q_{n-1}][t]$-linear operators
$\tilde{\sigma}_{ij}^q$ by 
\[ \tilde{\sigma}_{ij}^q (f(y)\underline{w})=t(\pa_{w(i)w(j)} f(y))\underline{w}+
f(y)\sigma^q_{ij}(\underline{w}) . \] 
We can check the well-definedness of the quantum extended Bruhat representation 
$[ij] \mapsto \tilde{\sigma}_{ij}^q,$ $x_k\mapsto \xi_k$ of the algebra $\E_n^q \langle R \rangle [t]$ 
in the same way as the proof of Proposition 3.1. 
\begin{thm} {\it 
Let $\S^q_w(x,y)$ be the quantum double Schubert polynomial corresponding to $w\in {\bf S}_n$ 
$($see {\rm \cite{CF}} and {\rm \cite{KM}).} 
When $t=0,$ we have
\[ \S^q_w(\theta,y)(\underline{\rm{id.}})=\underline{w} \] 
under the quantum extended Bruhat representation.}
\end{thm}
Proof.\quad This follows from the quantum Monk formula \cite{FGP} for the quantum Schubert
polynomials. 

\begin{cor} {\it 
The quantum double Schubert polynomials $\S^q_w(x,y)$ are characterized by the 
conditions: \\ 
$(1)$ $\S^q_w(x,y)|_{q=0}=\S_w(x,y),$ \\ 
$(2)$ $\S^q_w(x,y)$ is a linear combination of polynomials $\S_v(x,y)$ with $v \leq w$ over 
$\Z[q_1,\ldots , q_{n-1}],$ \\ 
$(3)$ $\S^q_w(\theta,y)(\underline{\rm{id.}})=\underline{w}$ under the quantum extended 
Bruhat representation at $t=0.$}  
\end{cor} 
Proof.\quad Let $\{ P_w(x,y) \}_{w \in {\bf S}_n}$ be a family of polynomials 
which satisfy the above properties $(1),$ $(2)$ and $(3)$ of the quantum 
Schubert polynomials. We will show that $P_w(x,y)=\S^q_w(x,y)$ for all 
$w\in {\bf S}_n.$ Let $S=1+(q_1,\ldots , q_{n-1})$ be a multiplicative system 
of $\Z[q_1,\ldots,q_{n-1}]$ consisting of polynomials of form 
$1+\sum_{i=1}^{n-1}q_i a_i(q),$ $a_i(q)\in \Z[q_1,\ldots,q_{n-1}].$ 
The properties $(1)$ and $(2)$ imply that 
\[ \S^q_w(x,y)= \sum_{v\leq w}b^v_w(q) \S_v(x,y), \;\;\; b^v_w(q)\in \Z[q_1,\ldots,q_{n-1}], \; 
b^w_w(q) \in S. \] 
Then it is easy to see that 
\[ \S_w(x,y)=\sum_{v\leq w}c^v_w(q) \S^q_v(x,y), \;\;\; 
c^v_w(q)\in S^{-1}\Z[q_1,\ldots,q_{n-1}], \] 
by induction on the Bruhat ordering. Similarly we also have 
\[ \S_w(x,y)=\sum_{v\leq w}d^v_w(q) P_v(x,y), \;\;\; 
d^v_w(q)\in S^{-1}\Z[q_1,\ldots,q_{n-1}]. \] 
Hence we have 
\[ \sum_{v\leq w}(c^v_w(q) \S^q_v(x,y)-d^v_w(q) P_v(x,y))=0 \] 
for all $w\in {\bf S}_n.$ From the property $(3),$ we obtain 
\[ \sum_{v\leq w}(c^v_w(q)-d^v_w(q)) \underline{v}=0 \] 
for all $w\in {\bf S}_n,$ so $c^v_w(q)=d^v_w(q)$ for all $w\in {\bf S}_n$ 
and $v\leq w.$ We can conclude that $P_w(x,y)=\S^q_w(x,y)$ from these identities 
again by induction on the Bruhat ordering. 

\section{Nichols-Woronowicz model}
The model of the equivariant cohomology ring $H^*_T(Fl_n)$ in the algebra
$\te_n$ has a natural interpretation in terms of the Nichols-Woronowicz
algebra. The Nichols-Woronowicz approach leads us to the uniform construction
for arbitrary root systems. For the definition of the Nichols-Woronowicz algebra, 
see e.g. \cite{AS}, \cite{Ba}, \cite{Maj} and \cite{Wo}. 

We denote by $\B_W$ the Nichols-Woronowicz algebra associated to
the Yetter-Drinfeld module
\[ V=\bigoplus_{\alpha \in \Delta}\R [\alpha]/([\alpha]+[-\alpha]) \]
over the finite Coxeter group $W$ of the root system $\Delta.$
Let $\h$ be the reflection representation of $W$ and $R={\rm Sym}\h^*$
the ring of polynomial functions on $\h.$
Let us consider the extension $\B_W\langle R \rangle [t]$ of the algebra
$\B_W$ by the
polynomial ring $R[t]$ defined by the commutation relation
\[ [\alpha]x= s_{\alpha}(x)[\alpha]+t(x,\alpha) \;\;\;\; \textrm{for} \; \;
x\in \h^*. \]
\begin{dfn} {\it 
We define the $R$-algebra $\tb_W$ by
\[ \tb_W= \B_W\langle R \rangle [t] \otimes_{R^{W}} R . \]}
\end{dfn}
Choose a $W$-invariant constants $(c_{\alpha})_{\alpha}.$
Let us consider a linear map $\mu : \h^* \rightarrow \tb_W$ defined
as
\[ \mu(x)= x+ \sum_{\alpha \in \Delta_+}c_{\alpha}(x,\alpha)[\alpha]  \]
for $x \in \h^*.$
\begin{prop}
$[\mu(x),\mu(y)]=0,$ $x,y\in \h^*.$
\end{prop}
Proof.\quad Let $\mu_0(x):=\sum_{\alpha \in \Delta_+}c_{\alpha}(x,\alpha)[\alpha].$ 
Here we may normalize the length of the roots to have $(\alpha,\alpha)=1,$ 
$\alpha \in \Delta.$ 
The commutativity $[ \mu_0(x), \mu_0(y)]=0$ has been shown in \cite{Ba}. 
The commutativity between $\mu(x)$ and $\mu(y)$ follows from 
\begin{eqnarray*}
\mu_0(x)y+x\mu_0(y) &=& \sum_{\alpha \in \Delta_+} c_{\alpha}(x,\alpha)s_{\alpha}(y)[\alpha] 
+ t \sum_{\alpha \in \Delta_+} c_{\alpha}(x,\alpha)(y,\alpha)
+ \sum_{\alpha \in \Delta_+} c_{\alpha}(y,\alpha)x [\alpha] \\ 
&=& \sum_{\alpha \in \Delta_+} c_{\alpha}(x,\alpha)y[\alpha] - \sum_{\alpha \in \Delta_+} 
2c_{\alpha}(x,\alpha)(y,\alpha)\alpha [\alpha] \\ 
& & + t \sum_{\alpha \in \Delta_+} c_{\alpha}(x,\alpha)(y,\alpha)
+ \sum_{\alpha \in \Delta_+} c_{\alpha}(y,\alpha)x [\alpha] \\ 
&=& \sum_{\alpha \in \Delta_+} c_{\alpha}(x,\alpha)y[\alpha]
+ t \sum_{\alpha \in \Delta_+} c_{\alpha}(x,\alpha)(y,\alpha)
+ \sum_{\alpha \in \Delta_+} c_{\alpha}(y,\alpha)s_{\alpha}(x) [\alpha] \\ 
&=& y\mu_0(x)+\mu_0(y)x.
\end{eqnarray*}

Proposition 5.1 shows that the linear map $\mu$ extends to a homomorphism of algebras
\[ \mu : R \rightarrow \B_W\langle R \rangle [t]. \]
Denote by $\widetilde{\mu}$ the composite of the homomorphisms 
\[ R \otimes_{\R} R \stackrel{\mu \otimes 1}{\rightarrow} \B_W\langle R \rangle [t] \otimes _{\R} R
\rightarrow \tb_W . \]

We will show in Theorem 5.1 below that the image of the algebra homomorphism $\widetilde{\mu}$ 
at $t=0$ is isomorphic the algebra $R\otimes_{R^W}R.$ 
The proof is based on the correspondence between the twisted
derivation $D_{\alpha}$ and the divided difference operator 
$\pa_{\alpha}:=(1-s_{\alpha})/\alpha,$ which acts on the first tensor component 
of $R\otimes_{R^W}R$ and extends linearly with respect to the second tensor 
component. 
We define the operator $D_{\alpha}$ as the twisted derivation on $\tb_W$
determined by the conditions: \\
(1): $D_{\alpha}(x)=0,$ for $x\in R,$ \\
(2): $D_{\alpha}([\beta])= \delta_{\alpha,\beta},$ for $\alpha,\beta \in \Delta_+,$ \\
(3): $D_{\alpha}(fg)=D_{\alpha}(f)g+s_{\alpha}(f)D_{\alpha}(g).$ \\ 
The operator $D_{\alpha}$ is linear with respect to $R$ on the second component. 
\begin{prop}
\[ \cap_{\alpha \in \Delta_+}{\rm Ker}(D_{\alpha})=R[t]\otimes_{R^W}R. \]
\end{prop}
Proof.\quad Since $\B_W\langle R \rangle [t] \cong R[t] \otimes_{\bf R} \B_W$ 
as a rignt $\B_W$-module, 
any element $\omega \in \B_W\langle R \rangle [t]$ can be written as 
\[ \omega = f_1 \varphi_1+ \cdots +f_k \varphi_k, \] 
where $f_1,\ldots,f_k \in R[t]$ are linearly independent, and 
$\varphi_1,\ldots ,\varphi_k \in \B_W.$ We have 
\[ D_{\alpha}(\omega)= s_{\alpha}(f_1)D_{\alpha}(\varphi_1)+\cdots + 
s_{\alpha}(f_k)D_{\alpha}(\varphi_k) \] 
from the twisted Leibniz rule. If $D_{\alpha}(\omega)=0,$ we have 
$D_{\alpha}(\varphi_1)=\cdots =D_{\alpha}(\varphi_k)=0.$ 
Hence, $\omega \in \cap_{\alpha \in \Delta_+}{\rm Ker}(D_{\alpha})$ 
implies that $\varphi_i$ belongs to the homogeneous part 
$\B^0_W \cong \R$ of degree zero for $i=1,\ldots ,k.$ 
This means $\omega \in R[t].$ 
\begin{prop} {\it 
\[ D_{\alpha}(\widetilde{\mu}(x))=c_{\alpha}\widetilde{\mu}(\partial_{\alpha}(x)) \]
for $x\in R \otimes_{\R} R.$}
\end{prop}
Proof.\quad When $x=\beta \otimes 1,$ $\beta \in \Delta,$ we can check 
that 
\[ D_{\alpha}(\widetilde{\mu}(\beta \otimes 1))= c_{\alpha}(\beta,\alpha) = 
c_{\alpha}\widetilde{\mu}(\partial_{\alpha}(\beta)) . \] 
Hence, we have $D_{\alpha}(\widetilde{\mu}(x))=
c_{\alpha}\widetilde{\mu}(\partial_{\alpha}(x))$ for $x\in \h^* \otimes R.$ 
On the other hand, the both-hands sides satisfy the same twisted 
Leibniz rule, so it follows that $D_{\alpha}(\widetilde{\mu}(x))=
c_{\alpha}\widetilde{\mu}(\partial_{\alpha}(x))$ for $x\in R \otimes R.$ 
\begin{thm} {\it 
If $t=0$ and the constants $(c_{\alpha})_{\alpha}$ are generic,
the image of the homomorphism $\widetilde{\mu}$ is isomorphic to 
the algebra $R \otimes_{R^W} R.$ In particular, when $W$ is the Weyl 
group, it is isomorphic to 
the $T$-equivariant cohomology ring $H^*_T(G/B)$ of the corresponding 
flag variety $G/B.$} 
\end{thm}
Proof.\quad If $x\in R^W \otimes_{\R} R,$ we have 
$D_{\alpha}(\widetilde{\mu}(x))=0$ for every 
$\alpha \in \Delta_+$ from Proposition 5.3. This implies from Proposition 5.2 
that $\widetilde{\mu}(x)\in R^W \otimes_{R^W} R.$ When $t=0,$ 
$\widetilde{\mu}(x)$ coincides with the element of $R$ which is 
obtained by replacing all the symbols $[\alpha]$ by zero in 
$\widetilde{\mu}(x).$ Hence, 
the homomorphism $\widetilde{\mu}$ factors through $R\otimes_{R^W}R 
\rightarrow \tb_W.$ 

Let us take a reduced decomposition $w=s_{\alpha_{i_1}}\cdots s_{\alpha_{i_l}}$ 
of an element $w\in {\bf S}_n.$ Then the operator 
$\partial_w:= \partial_{\alpha_{i_1}}\cdots \partial_{\alpha_{i_l}}$ is 
independent of the choice of reduced decompositions thanks to the Coxeter 
relation. We also define $D_w:=D_{\alpha_{i_1}}\cdots D_{\alpha_{i_l}}.$ 
Define a family $\{ X_w \}_{w\in {\bf S}_n}$ of polynomials 
by $X_{w_0}:=|W|^{-1}\prod_{\alpha \in \Delta_{+}}\alpha$ and 
$X_w:=\partial_{w^{-1}w_0}X_{w_0}.$ The family $\{ X_w \}_{w\in W}$ gives 
a linear basis of $R_W.$ We can see that 
\[ CT(\partial_w X_v)= \left\{ 
\begin{array}{cc} 
1 & \textrm{if $w=v,$} \\ 
0 & \textrm{otherwise.} 
\end{array}
\right. \] 
where $CT$ stands for the part of degree zero. 
Since a linear basis $\{ X_w \}_{w\in W}$ of the coinvariant algebra of 
$W$ gives an $R^W$-basis of $R,$ and $D_w\widetilde{\mu}(X_v \otimes 1)= 
c_{\alpha_{i_1}}\cdots c_{\alpha_{i_l}}\widetilde{\mu}(\partial_w X_v \otimes 1),$ 
it is easy to see that $R\otimes_{R^W}R 
\rightarrow \tb_W$ is injective.

\begin{rem}
Our construction is not a straightforward application of the functor 
$(-)\otimes_{R^W}R$ to the one given in \cite{Ba} even when $t=0.$ 
In fact, the defining relations of the algebra $\B_W\langle R \rangle [t]$ 
involve a nontrivial commutation relation 
\[ [\alpha]x= s_{\alpha}(x)[\alpha]+t(x,\alpha). \] 
Moreover, in the formula 
\[ \mu(x)= x+ \sum_{\alpha \in \Delta_+}c_{\alpha}(x,\alpha)[\alpha], \] 
the first term in the right-hand side does not appear in the non-equivariant case. 
\end{rem}

Anatol N. Kirillov \\
Research Institute for Mathematical Sciences \\ 
Kyoto University \\
Sakyo-ku, Kyoto 606-8502, Japan \\
e-mail: {\tt kirillov@kurims.kyoto-u.ac.jp} \\
URL: {\tt http://www.kurims.kyoto-u.ac.jp/\textasciitilde kirillov} 
\bigskip 
\\ 
Toshiaki Maeno \\ 
Department of Electrical Engineering \\ 
Kyoto University \\ 
Sakyo-ku, Kyoto 606-8501, Japan \\ 
e-mail: {\tt maeno@kuee.kyoto-u.ac.jp}
\end{document}